\documentclass[12pt]{article}
\usepackage{amsmath}
\usepackage{amssymb}
\usepackage{amscd}
\usepackage{amsthm}

\theoremstyle{plain}
\newtheorem{Thm}{Theorem}[section]

\newtheorem{Prop}[Thm]{Proposition}
\newtheorem{Cor}[Thm]{Corollary}
\newtheorem{Lem}[Thm]{Lemma}

\theoremstyle{definition}
\newtheorem{Defn}[Thm]{Definition}
\newtheorem{Expl}[Thm]{Example}

\newtheorem{Que}[Thm]{Question}

\numberwithin{equation}{section}

\title{Francia's flip and derived categories}
\author{Yujiro Kawamata}

\begin{document}

\maketitle

\begin{abstract}
We extend some of the results of Bondal-Orlov on the equivalence of 
derived categories to the case of orbifolds 
by using the category of coherent orbifold sheaves.
\end{abstract}

\section{Introduction}

We consider an approach to the problems on flips and flops in the birational 
geometry from the point of view of the theory of derived categories. 
The purpose of this paper is to extend some of the existing results for smooth 
varieties to the case of varieties having only quotient singlarities.

The idea of using the derived categories can be explained in the following way.
The category of sheaves on a given variety is directly related to the 
biregular geometry of the variety.  
But the category derived from the category of complexes of sheaves 
by adding the inverses of quasi-isomorphisms and dividing modulo chain
homotopy equivalences acquires more symmetry, and is beleived to reflect more 
essential properties of the variety, namely the birational geometry of 
the variety.
More precisely, the varieties which have the same level of the 
canonical divisors
(so called $K$-equivalent varieties) are beleived to have equivalent 
derived categories.

Bondal-Orlov \cite{BO1} considered a smooth variety $X$ which contains a 
subvariety $E$ isomorphic to the projective space $\mathbb{P}^m$ 
such that the normal bundle is isomorphic to $\mathcal{O}(1)^{n+1}$.
If we blow up $X$ with center $E$, then the exceptional divisor can be 
contracted 
to another direction to yield another smooth variety $X'$ which contains a 
subvariety $E'$ isomorphic to the projective space $\mathbb{P}^n$.
The induced birational map $X -\to X'$ is a flip if $m > n$ and a flop 
if $m=n$.
In other words, the $K$-level of $X$ is higher than that of $X'$ if $m>n$ and 
they are equal if $m=n$.
Then \cite{BO1} proved that the natural functor 
between the derived categories of
bounded complexes of coherent sheaves 
$D^b_{\text{coh}}(X') \to D^b_{\text{coh}}(X)$, called the Fourier-Mukai
transform after \cite{Mukai},
is fully faithful if $m > n$ and an equivalence if $m=n$.

\cite{BO1} (see also \cite{BO2}) also proved
a reconstruction theorem in the following sense:
if there exists an equivalence of derived categories
$D^b_{\text{coh}}(X) \to D^b_{\text{coh}}(X')$
for smooth projective varieties $X$ and $X'$ such that either the canonical 
divisor $K_X$ or its negative $-K_X$ is ample, 
then there exists an isomorphism $X \to X'$.

We shall extend these results for varieties having quotient singularities in
this paper.
The main results are Theorem~\ref{main1} and Theorem~\ref{main2}.
We consider some toric flips and flops defined in \S4 after \cite{Reid} and 
\cite{Thaddeus}.
We first remark in Example~\ref{counterexample} of \S5 
that this kind of extension
does not work if we consider the usual derived categories of 
bounded complexes of coherent sheaves.
In order to overcome this difficulty, we introduce the concept of coherent 
orbifold sheaves in \S2. 
In Theorem~\ref{main1}, we prove that only the level of $K$ deternimes
the equivalence class of the derived categories, though the varieties 
with the same level of $K$ may have very different geometric outlook.
For example, the dimensions of the exceptional loci may be different.
Unlike the smooth case, there is no obvious geometric order 
between the varieties, though there is order of canonical 
divisors, and the derived categories follow the latter.

According to the minimal model program, we should deal with varieties which 
admit mild singularities, 
and results for smooth varieties should be extended to such varieties 
(cf. \cite{KMM}).
Our extension for varieties with quotient singularities can be regarded 
as the first step toward the general case of varieties with arbitrary terminal 
singularities.
A recent result by Yasuda \cite{Yasuda} on the motivic integration for 
orbifolds is also one of such extensions.

The existence of the flips for arbitrary small contraction with relatively 
negative canonical divisor is one of the most important but difficult 
conjectures in the minimal model program.
It is proved only in dimension $3$ by Mori \cite{Mori}.
In \S3, we recall a result of the author \cite{Kawamata} which reduces 
the existence problem of the flips to that of the flops 
(Theorem~\ref{flipflop}).
The reason is that the flops seem to be better suited to the categorical 
argument than the flips,
because the flop corresponds to the equivalence of categories while the
flip to the fully faithful embedding. 

Bridgeland \cite{Bridgeland2} constructed the flop for any small crepant 
contraction of a smooth $3$-dimensional variety
by using only the categorical argument as in \cite{BKR}.
It is remarkable that the existence of the flop and the equivalence of 
derived categories are simultaneously proved.
While preparing this manuscript, the author learned that Chen \cite{C} 
announced a result which extends the above result \cite{Bridgeland2}
to the flops of $3$-dimensional varieties with Gorenstein terminal 
singularities.
One might even extend this to the case of $3$-folds having arbitrary 
terminal singularities by combining 
with our method since such singularities can be deformed to quotient 
singularities.  
We hope that we could eventually prove the existence of flips in this way. 

The author would like to thank Akira Ishii and Adrian Langer 
for the useful discussions on the
derived categories and orbifold sheaves, respectively.
We work over the complex number field $\mathbb{C}$.


\section{Orbifold sheaf}

We begin with recalling the definition of the quasi-projective 
orbifolds (or $Q$-varieties) and coherent orbifold sheaves (or $Q$-sheaves)
from \cite{Mumford}~\S2.

\begin{Defn}
Let $X$ be a quasi-projective variety.
An {\it orbifold structure} on $X$ consisits of the 
data $\{\pi_i: X_i \to X, G_i\}_{i \in I}$,
where the $X_i$ are smooth quasi-projective varieties,
the $\pi_i$ are quasi-finite morphisms, and the $G_i$ are finite 
groups acting faithfully on the $X_i$, such
that $X = \bigcup_{i \in I} \pi_i(X_i)$, the $\pi_i$ induce etale morphisms 
$\pi_i': X_i/G_i \to X$, and that, if $X_{ij} = (X_i \times_X X_j)^{\nu}$ 
denotes the normalization of the fiber product, then 
the projections $p_1: X_{ij} \to X_i$ and $p_2: X_{ij} \to X_j$ are 
etale for any $i$ and $j$, where $i$ and $j$ may be equal. 
\end{Defn}

In this case, $X$ has only quotient singularities.  Conversely, 
if $X$ is a quasi-projective variety having only quotient singularities, 
then there exists an orbifold structure on $X$ such that the $p_i$ 
are etale in codimension $1$. We call such a structure {\it natural}.

A {\it global cover} $\tilde X$ is the normalization of $X$ in a 
Galois extension of
the function field $k(X)$ which contains all the extensions $k(X_i)$.

\begin{Defn}
An {\it orbifold sheaf} $F$ is a collection of sheaves $F_i$
of $\mathcal{O}_{X_i}$-modules on the $X_i$ together 
with the gluing isomorphisms 
$g_{ji}: p_1^*F_i \to p_2^*F_j$ on the $X_{ij}$, 
such that the compatibility conditions 
$(p_{23}^*g_{kj}) \circ (p_{12}^*g_{ji}) = p_{13}^*g_{ki}$ 
hold on the triple overlaps
$X_{ijk} = (X_i \times_X X_j \times_X X_k)^{\nu}$, where $\nu$ denotes
the normalization.
\end{Defn}

For example, we define the {\it orbifold structure sheaf}
$\mathcal{O}_X^{\text{orb}}$ by the $\mathcal{O}_{X_i}$.

Let $\tilde X_i$ be the normalization of $X_i$ in the function field 
$k(\tilde X)$ of the global cover, and 
$H'_i = \text{Gal}(\tilde X_i/X'_i)$, where $X'_i = X_i/G_i$.
Then an orbifold sheaf $F$ on $X$ is in a one-to-one correspondense 
to a sheaf $\tilde F$ of $\mathcal{O}_{\tilde X}$-modules on $\tilde X$ 
such that the action of the Galois group $\text{Gal}(\tilde X/X)$ lifts 
to $\tilde F$ and that the restriction $\tilde F \vert_{\tilde X_i}$ 
with its $H_i'$-action is isomorphic to the pull-back of 
a sheaf of $\mathcal{O}_{X_i}$-modules on $X_i$.  

A {\em homomorphism} $h: F \to F'$ of orbifold sheaves is a collection of 
$G_i$-equivariant $\mathcal{O}_{X_i}$-homomorphisms $h_i: F_i \to F'_i$ 
which are compatible with gluing isomorphisms.
A {\it tensor product} $F \otimes_{\mathcal{O}_X^{\text{orb}}} F'$ is 
given by the sheaves $F_i \otimes_{\mathcal{O}_{X_i}} F'_i$.
The category of orbifold sheaves $Sh(X^{\text{orb}})$ on $X$ 
thus defined becomes an abelian category.

An orbifold sheaf is said to be {\it coherent} (resp. {\it locally free}) 
if each $F_i$ is coherent (resp. locally free).
For example, the {\it orbifold sheaf of differential $p$-forms} 
$\Omega_X^{p,\text{orb}}$ is a locally free coherent orbifold sheaf 
given by the sheaves $\Omega_{X_i}^p$.
In particular, the {\it dualizing orbifold sheaf} $\omega_X^{\text{orb}}$ 
is the invertible orbifold sheaf consisting of the dualizing 
sheaves $\omega_{X_i}$. 
We note that even if $\omega_{X_i}$ is isomorphic to $\mathcal{O}_{X_i}$, 
the action of $G_i$ on them may be different.

Let $D^b(X_{\text{coh}}^{\text{orb}})$ (resp. 
$D^b_c(X_{\text{coh}}^{\text{orb}})$) be the derived category of 
bounded complexes of coherent orbifold sheaves on $X$
(resp. with compact supports).

\begin{Prop}\label{resolutions}
Let $X$ be a quasi-projective variety with an orbifold structure.

(1) The category $Sh(X^{\text{orb}})$ has enough injectives.

(2) (\cite{Mumford}~Proposition~2.1.) 
If $X$ has a Cohen-Macaulay global cover,
then any coherent orbifold sheaf $F$ has a finite locally free resolution.
\end{Prop}

\begin{Prop}\label{images}
Let $f: X \to Y$ be a generically surjective 
morphism of quasi-projective varieties with orbifold structures
$\{\pi_i: X_i \to X, G_i\}$ and 
$\{\rho_{\alpha}: Y_{\alpha} \to Y, H_{\alpha}\}$.
Assume that the natural morphism 
$p_1^{i\alpha}: (X_i \times_Y Y_{\alpha})^{\nu} \to X_i$ 
is etale for any $i$ and $\alpha$, where $\nu$ denotes the normalization. 
Then one can define the direct image functor 
$f_*: Sh(X^{\text{orb}}) \to Sh(Y^{\text{orb}})$ and
the inverse image functor
$f^*: Sh(Y^{\text{orb}}) \to Sh(X^{\text{orb}})$ 
which are adjoints each other.
\end{Prop}

\begin{proof}
When we fix $\alpha$ and vary $i$, then we obtain a covering of the normalized 
fiber product $(X \times_Y Y_{\alpha})^{\nu}$ by the morphisms
$\sigma_{i\alpha}: (X_i \times_Y Y_{\alpha})^{\nu} \to 
(X \times_Y Y_{\alpha})^{\nu}$ induced from the $\pi_i$.
For an orbifold sheaf $E$ on $X$, we define a sheaf $E_{\alpha}$ on 
$(X \times_Y Y_{\alpha})^{\nu}$ as the kernel
\[
E_{\alpha} \to \bigoplus_i \sigma_{i\alpha *}p_1^{i\alpha*}E_i 
\rightrightarrows 
\bigoplus_{i,j} \sigma_{ij\alpha *}p_1^{ij\alpha*}E_i
\]
where $p_1^{ij\alpha}: (X_i \times_X X_j \times_Y Y_{\alpha})^{\nu} \to X_i$ 
and $\sigma_{ij\alpha}: (X_i \times_X X_j \times_Y Y_{\alpha})^{\nu} \to
(X \times_Y Y_{\alpha})^{\nu}$ are natral morphisms.
Then we define an orbifold sheaf $f_*E$ on $Y$ by 
$(f_*E)_{\alpha} = p_{2*}E_{\alpha}$.
Since the $p_1^{i\alpha}$ are etale, $f_*$ is left exact. 

When we fix $i$ and vary $\alpha$, then we obtain an etale covering  
of $X_i$ by $p_1^{i\alpha}: (X_i \times_Y Y_{\alpha})^{\nu} \to X_i$.
For an orbifold sheaf $F$ on $Y$, we define a sheaf $(f^*F)_i$ 
on $X_i$ as the kernel
\[
(f^*F)_i \to \bigoplus_{\alpha} p_{1*}^{i\alpha}p_2^{i\alpha*}F_{\alpha} 
\rightrightarrows \bigoplus_{\alpha, \beta} p_{1*}^{i\alpha \beta}
p_2^{i\alpha \beta*}F_{\alpha}
\]
where $p_1^{i\alpha \beta}: (X_i \times_Y Y_{\alpha} 
\times_Y Y_{\beta})^{\nu} \to X_i$ and
$p_2^{i\alpha \beta}: (X_i \times_Y Y_{\alpha} 
\times_Y Y_{\beta})^{\nu} \to Y_{\alpha}$ are natural morphisms.
Since the $p_1^{i\alpha}$ are etale, $f^*$ is right exact.
\end{proof}

\begin{Cor}\label{derived-image}
In addition to the assumptions of Proposition~\ref{images},
assume that $f$ is proper and $Y$ has a Cohen-Macaulay global cover.
Then functors $f_*$ and $f^*$ induces exact functors
\[
\begin{split}
&Rf_*^{\text{orb}}: D^b(X_{\text{coh}}^{\text{orb}}) 
\to D^b(Y_{\text{coh}}^{\text{orb}}) \\
&Lf^*_{\text{orb}}: D^b(Y_{\text{coh}}^{\text{orb}}) 
\to D^b(X_{\text{coh}}^{\text{orb}}).
\end{split}
\]
\end{Cor}

We reark that we need to consider $Q$-stacks instead of $Q$-varieties
in order to deal with general morphisms $f: X \to Y$ of orbifolds.
 
We have the Serre functor as follows:

\begin{Prop}
Let $X$ be a quasi-projective variety with an orbifold structure.
Assume that $X$ has a Cohen-Macaulay global cover $\tilde X$.
Then there exists a Serre functor $S = S_X$ for the derived categories
$D = D^b(X_{\text{coh}}^{\text{orb}})$ and 
$D_c = D^b_c(X_{\text{coh}}^{\text{orb}})$
defined by 
\[
S(u) = u \otimes_{\mathcal{O}_X^{\text{orb}}} \omega_X^{\text{orb}}[\dim X]
\]
for $u \in D$.
There are bifunctorial isomorphisms
\[
\text{Hom}_D(u, v) 
\cong \text{Hom}_D(v, S(u))^*
\]
for $u \in D$ and $v \in D_c$ or $u \in D_c$ and $v \in D$.
\end{Prop}

\begin{proof}
We first assume that $u$ is a locally free orbifold sheaf and $v$ is an orbifold sheaf 
with compact support.
If we replace $v$ by $u^* \otimes v$, we may assume that
$u = \mathcal{O}_X^{\text{orb}}$.
We denote by $\tilde v$ the sheaf on $\tilde X$ corresponding to $v$
with the action of $G = \text{Gal}(\tilde X/X)$.
Let $\pi: \tilde X \to X$ denote the natural morphism.
We have 
\[
\text{Hom}_D(u,v[k])
\cong H^k(\tilde X, \tilde v)^G 
\cong H^k(X, (\pi_*\tilde v)^G). 
\]
On the other hand, the Zariski sheaves $v'_i = (\pi_{i*}v_i)^{G_i}$ on the 
$X'_i = X_i/G_i$ define an etale sheaf on $X$. Indeed, we have
$v'_i = (\pi'_i)^*(\pi_*\tilde v)^G$, where $\pi'_i: X'_i \to X$ is induced from $\pi_i$.
By the relative duality for the morphism $X_i \to X'_i$, we have
$RHom(v_i, \omega_{X_i})^{G_i} \cong RHom(v'_i, \omega_{X'_i})$.
Since $\omega_{X'_i} = (\pi'_i)^*\omega_X$, we have 
\[
\text{Hom}_D(v[k], S(u))
\cong \text{Ext}^{d-k}((\pi_*\tilde v)^G, \omega_X)
\]
where $d = \dim X$.
Hence our assertion is reduced to the usual duality theorem on $X$.

If $v$ is a locally free orbifold sheaf and $u$ is an orbifold sheaf 
with compact support, then by the first part
\[
\text{Hom}_D(u,v[k])
\cong \text{Hom}_D(S(u),S(v)[k])
\cong \text{Hom}_D(v[k],S(u))^*.
\]
The general case is obtained by taking the locally free resolutions.
\end{proof}


\section{Flip to flop}

We shall reduce the existence problem of the flips to that of flops.
For this purpose, 
we consider the total space $KX$ of the $\mathbb{Q}$-bundle $K_X$:

\begin{Lem}\label{cone}
Let $X$ be a variety of dimension $n$ 
with only log terminal singularities.
Then 
\[
KX = \text{Spec } (\bigoplus_{m=0}^{\infty} \mathcal{O}_X(- mK_X)) 
\]
is a variety of dimension $n+1$ 
with only rational Gorenstein singularities and trivial canonical bundle.
\end{Lem}

\begin{proof}
We take a small open subset $U$ of $X$ which has an index $1$ cover
$\pi: U_1 \to U$ with the Galois group $G$.
We have a commutative diagram
\[
\begin{CD}
KU_1 @>{K\pi}>> KU \\
@VVV          @VVV \\
U_1 @>{\pi}>> U.
\end{CD}
\]
$U_1$ has only rational Gorenstein singularities and $KU_1$ is the total 
space of the line bundle $K_{U_1}$, hence $KU_1$ has only 
rational Gorenstein singularities.
Since $\mathcal{O}_X(- mK_U) = (\pi_*\mathcal{O}_{U_1}(- mK_{U_1}))^G$,
we have $KU_1 = KU/G$.

Let $\omega$ be a generating section of $K_{U_1}$.
Then $\omega^{-1}$ gives a fiber coordinate along the fiber of $K_{U_1}$,
and $\tilde{\omega} = \omega \wedge d \omega^{-1}$ is a 
generating section of $K_{KU_1}$.
Since $\tilde{\omega}$ is $G$-invariat, $K_{KX}$ is again invertible. 

Let $U'$ be another small open subset of $X$ 
which has an index $1$ cover $\pi': U'_1 \to U'$ and $\omega'$ a 
generating section of $K_{U'_1}$.
We can write $\omega' = u\omega$ for an invertible function $u$ on
$U_1 \times_X U'_1$.
Then $\tilde{\omega'} = \omega' \wedge d\omega^{\prime -1}
= u\omega \wedge u^{-1}d \omega^{-1} = \omega \wedge d \omega^{-1}
=\tilde{\omega}$.
Therefore, we obtain a global generating section of $K_{KX_1}$.
\end{proof}

\begin{Que}
If $X$ has only terminal singularities, so has $KX$?
\end{Que}

\begin{Thm}\label{flipflop}
Let $n$ be an integer such that $n \ge 3$.
Assume that the Gorenstein canonical flop always exists in dimension $n+1$ 
in the following sense:
for any projective birational small crepant morphism $\psi: \tilde X
\to \tilde Z$ of $(n+1)$-dimensional varieties with only Gorenstein canonical 
singularities and a $\psi$-negative effective $\mathbb{Q}$-Cartier divisor
$D$ on $\tilde X$, 
there exists another projective birational small crepant morphism 
$\psi^+: \tilde X^+ \to \tilde Z$ from a $(n+1)$-dimensional variety 
with only Gorenstein canonical singularities such that the strict transform
$D^+$ of $D$ on $\tilde X^+$ is $\psi^+$-ample.
Then the canonical flip always exists in dimension $n$ 
in the following sense:
for any projective birational small morphism $\phi: X
\to Z$ from an $n$-dimensional variety with only canonical 
singularities such that $K_X$ is $\phi$-negative, 
there exists another projective birational small morphism 
$\phi^+: X^+ \to Z$ from an $n$-dimensional variety 
with only canonical singularities such that $K_{X^+}$ is $\phi^+$-ample.
\end{Thm}

\begin{proof}
Let $\phi: X \to Z$ be a small projective birational morphism from 
a variety of dimension $n$ 
with only canonical singularities to a normal variety 
such that the canonical divisor $K_X$ is $\phi$-negative.
Let $E$ be the exceptional locus of $\phi$.
Let $KX$ be as in Lemma~\ref{cone} with a natural morphism
$\xi: KX \to X$.
We have an embedding $i: X \to KX$ to the zero section.
We also define
\[
KZ = \text{Spec } (\bigoplus_{m=0}^{\infty} \mathcal{O}_Z(- mK_Z)).
\]
Since 
\[
\mathcal{O}_Z(- mK_Z) \cong \phi_*(\mathcal{O}_X(- mK_X))
\] 
the direct sum on the right hand side of the formula for $KZ$ gives
a finitely generated sheaf of $\mathcal{O}_Z$-algebras.

We claim that the induced morphism $K\phi: KX \to KZ$ is a projective 
birational morphism whose exceptional locus coincides with $i(E)$.
In order to prove this, we may assume that $Z$ is affine.
In this case, $K\phi$ is given by the linear system $\Lambda$ 
generated by divisors
of the form $mi(X) + \xi^*D$ for $D \in \vert -mK_X \vert$.
Since $\vert -mK_X \vert$ is very ample on $X$ for sufficiently large $m$,
so is $\Lambda$ on $KX \setminus i(X)$.
Since the restriction of $K\phi$ to $i(X)$ coincides with $\phi$, we have our
assertion. 
Since $K_{KX}$ is globally trivial, 
$K\phi$ is crepant.

By the assupmtion, there exists a flop $(K\phi)^+: (KX)^+ \to KZ$
with respect to $i(X)$.
Let $X^+$ be the strict transform of $i(X)$ on $(KX)^+$.
By definition, $X^+$ is a $\mathbb{Q}$-Cartier divisor on $(KX)^+$ which is 
ample for $(K\phi)^+$.
Hence $X^+$ is $\mathbb{Q}$-Gorenstein though $X^+$ may not be normal.
We shall prove that $X^+$ is regular in codimension $1$ and the
induced morphism $\phi^+: X^+ \to Z$ is small.
Then $X^+$ is normal, since it is Cohen-Macaulay, 
and $K_{X^+} = K_{(KX)^+} + X^+ \vert_{X^+}$ is
$\phi^+$-ample, so that $\phi^+$ is the flip of $\phi$.

Let $\tilde Y$ be a common desingularization of $KX$ and $(KX)^+$ with 
projective birational morphisms $\tilde{\mu}: \tilde Y \to KX$ and
$\tilde{\mu}^+: \tilde Y \to (KX)^+$ with normal crossing exceptional locus
$F = \sum F_j$.
We write $\tilde{\mu}^*i(X) = Y + \sum_j r_jF_j$ and
$\tilde{\mu}^*K_{KX} + \sum_j a_jF_j= K_{\tilde Y}$, where 
$Y$ is the strict transform of $i(X)$.
Since $X$ has only canonical singularities,
we have $r_j \le a_j$ for any $j$.
We write also $(\tilde{\mu}^+)^*X^+ = Y + \sum_j r_j^+F_j$ and
$(\tilde{\mu}^+)^*K_{(KX)^+} + \sum_j a_j^+F_j= K_{\tilde Y}$.
Since $\tilde{\mu}^+$ is the flop of $\tilde{\mu}$ with respect to $i(X)$,
we have $r_j^+ < r_j$ for any divisor $F_j$ which lies above $i(E)$.
Since we have also $a_j^+ = a_j$ for any $j$,
we conclude that $r_j^+ < a_j^+$ for such $j$.

Let $P$ be any codimension $2$ point of $(KX)^+$ contained in 
the exceptional locus of $\tilde{\mu}^+$.
Since $(KX)^+$ has only canonical singularities,
there exists an exceptional divisor $F_{j_0}$ above $P$ such that
$a_{j_0}^+ = 0$ if $(KX)^+$ is singular at $P$ or $=1$ 
if $(KX)^+$ is smooth at $P$.
Thus $r_{j_0}^+ = 0$, and $X^+$ does not contain $P$.
Therefore, $X^+$ is regular in codimension $1$ and the
induced morphism $\phi^+: X^+ \to Z$ is small.
\end{proof}

\begin{Expl}
In the notation of \S4,
if $X^- = X^-(\mathbf a; \mathbf b)$ and 
$c = \sum_i a_i - \sum_j b_j > 0$, then
we have $KX^- = X^-(\mathbf a; \mathbf b,c)$.

Let $a$ and $b$ be coprime positive integers. 
The sequene of integers 
\[
(1,a,b;1,a+b)
\]
is obtained in the above way from the 
following sequences:
(1) $(1,a,b;a+b)$, (2) $(1,a,b;1)$, 
(3) $(1,a+b;a,b)$, (4) $(1,a+b;1,a)$, (5) $(1,a+b;1,b)$.
(1) and (2) correspond to divisorial contractions of $3$-folds with only
terminal quotient singularities to the quotient 
singularity of type $\frac 1{a+b}(1,a,b)$ and to a smooth point, respectively.
(3), (4) and (5) correspond to flips from $3$-folds with only
terminal quotient singularities.
\end{Expl}


\section{Toric flip and flop}

According to \cite{Reid} and \cite{Thaddeus}, we consider the following
toric varieties:

\begin{Defn}
Let $(\mathbf a; \mathbf b) = (a_1, \dots, a_m; b_1, \dots, b_n)$ be 
a sequence of positive integers.
We let the multiplicative group $G_m$ act on 
$\mathbb{A} = \mathbb{A}(\mathbf a; \mathbf b) = 
\text{Spec }R \cong \mathbb{A}^{m+n}$ 
for $R = \mathbb{C}[x_1,\dots,x_m,y_1,\dots,y_n]$ by
\[
\lambda_t(x_1, \dots, x_m, y_1, \dots, y_n) =  
(t^{a_1}x_1, \dots, t^{a_m}x_m, t^{-b_1}y_1, \dots, t^{-b_n}y_n)
\]
for $t \in G_m$. 
We consider GIT quotients 
\[
\begin{split}
&X^- = X^-(\mathbf a; \mathbf b) =
(\mathbb{A} \setminus \{x_1 =  \dots =  x_m = 0\})/G_m \\ 
&X^+ = X^+(\mathbf a; \mathbf b) = 
(\mathbb{A} \setminus \{y_1 =  \dots = y_n = 0\})/G_m \\
&X^0 = X^0(\mathbf a; \mathbf b) =\mathbb{A}//G_m = 
\text{Spec }R^{G_m}.
\end{split}
\] 
We also define $Y = Y(\mathbf a; \mathbf b)$ to be the fiber product 
$X^- \times_{X^0} X^+$.
Let $\phi^{\pm}: X^{\pm} \to X^0$ and $\mu^{\pm}: Y \to X^{\pm}$ 
be the induced morphisms as in the following commutative diagram:
\[
\begin{CD}
Y @>{\mu^+}>> X^+ \\
@V{\mu^-}VV   @VV{\phi^+}V \\
X^- @>>{\phi^-}> X^0
\end{CD}
\]
\end{Defn}

Let $A_i^{\pm}$ and $B_j^{\pm}$ be prime divisors on $X^{\pm}$ 
corresponding to the $x_i$ and the $y_j$,
and let $U_i^- = X^- \setminus A_i^-$ and $U_j^+ = X^+ \setminus B_j^+$.
Thus $X^- = \bigcup_i U_i^-$ and $X^+ = \bigcup_j U_j^+$.
Let $A_i$ and $B_j$ be the strict transforms of $A_i^{\pm}$ and $B_j^{\pm}$ 
on $Y$, respectively.
Let $U_{i,j} = Y \setminus (A_i \cup B_j)$ so that 
$Y = \bigcup_{i,j} U_{i,j}$.

\begin{Expl}
If $n = 0$, then $X^-$ is nothing but the weighted projective space
$\mathbb{P}(\mathbf a)$. In this case,
$X^0$ is a point and $X^+ = Y = \emptyset$.
If $n=1$, then $\phi^-$ is a divisorial contraction and $X^+ = X^0$.
If $m=n=2$ and $a_1=a_2=b_1=b_2=1$, then this is Atiyah's flop.
If $m=n=2$ and $a_1=2$, $a_2=b_1=b_2=1$, then this is Francia's flip.
\end{Expl}

\begin{Prop}
If $m,n \ge 2$, then the following hold.

(1) The morphisms $\phi^{\mp}$ are projective and birational whose  
exceptional loci $E^{\mp}$ are isomorphic to the weighted projective spaces
$\mathbb{P}(a_1, \dots, a_m)$ and $\mathbb{P}(b_1, \dots, b_n)$,
respectively.

(2) $E = (\mu^{\pm})^{-1}(E^{\pm})$ is a prime divisor on $Y$ isomorphic 
to the product $E^- \times E^+$.

(3) The divisors $\mp A_i^{\pm}$ and $\pm B_j^{\pm}$ are $\phi^{\pm}$-ample.
\end{Prop}

We have the following reduction for the sequence of integers 
in a similar way to the case of the weighted projective spaces 
(\cite{Delorme} and \cite{Dolgachev}).

\begin{Prop}\label{coprime}
(1) Let $c = GCD(a_1,\dots,a_m,b_1,\dots,b_n)$ be the greatest common divisor,
and let $(\mathbf a'; \mathbf b') = 
(a_1/c, \dots, a_m/c; b_1/c, \dots, b_n/c)$.
Then $X^{\pm}(\mathbf a'; \mathbf b') \cong X^{\pm}(\mathbf a; \mathbf b)$ and 
$X^0(\mathbf a'; \mathbf b') \cong X^0(\mathbf a; \mathbf b)$.

(2) Assume $c = 1$.
Let 
\[
\begin{split}
&c_i = GCD(a_1,\dots,\widehat{a_i}, \dots,a_m,b_1,\dots,b_n) \\
&c_{m+j} = GCD(a_1,\dots,a_m,b_1,\dots,\widehat{b_j}, \dots,b_n) \\ 
&d_k = LCM(c_1,\dots,\widehat{c_k}, \dots,c_{m+n}) \\
&(\mathbf a''; \mathbf b'') = (a_1/d_1, \dots, a_m/d_m; b_1/d_{m+1}, 
\dots, b_n/d_{m+n}).
\end{split}
\]
Then $X^{\pm}(\mathbf a''; \mathbf b'') \cong X^{\pm}(\mathbf a; \mathbf b)$ 
and 
$X^0(\mathbf a''; \mathbf b'') \cong X^0(\mathbf a; \mathbf b)$.
\end{Prop}

We may therefore assume that $c_i = c_{m+j} = 1$ for any 
$i$ and $j$ from now on.

\begin{proof}
(1) The action of $G_m$ factors through a homomorphism
$G_m \to G_m$ given by $t \mapsto t^c$. 

(2) Let $d = LCM(c_1,\dots,c_{m+n})$.
Since $c = 1$, we have $GCD(c_i,d_i) = 1$ and $d = c_id_i$.
The ring of invariants $R^{G_m}$ for the sequence $(\mathbf a; \mathbf b)$ is 
generated by monomials $x^{\mathbf m}y^{\mathbf n}$ such that 
$\sum_i a_im_i = \sum_j b_jn_j$, and
the corresponding ring of invariants $(R'')^{G_m}$ for the new sequence 
$(\mathbf a''; \mathbf b'')$ is generated by 
monomials $(x'')^{\mathbf m''}(y'')^{\mathbf n''}$ such that 
$\sum_i a_ic_im''_i = \sum_j b_jc_{m+j}n''_j$. 
If $\sum_i a_im_i = \sum_j b_jn_j$, then it follows that
$c_i \vert m_i$ and $c_{m+j} \vert n_j$.
Hence there is an isomorphism $f: (R'')^{G_m} \to R^{G_m}$ 
given by $f(x''_i) = x_i^{c_i}$ and $f(y_j'') = y_j^{c_{m+j}}$.
Thus we have $X^0(\mathbf a''; \mathbf b'') \cong X^0(\mathbf a; \mathbf b)$.
The assertions for $X^{\pm}$ follow from this isomorphism.
\end{proof}

\begin{Prop}
The open subset $U_1^-$ is isomorphic to the quotient of the 
affine space $\mathbb{A}^{m+n-1}$
by the group $\mathbb{Z}_{a_1}$ whose action is given by the weights
\[
\frac 1{a_1}(- a_2, \dots, -a_m, b_1, \dots, b_n).
\]
Moreover, the action of the group is small in the sense that 
the induced morphism $\mathbb{A}^{m+n-1} \to U_1^-$ is etale in codimension 
$1$.
\end{Prop}

\begin{proof}
Let $(\mathbf 1, \mathbf 1) = (1, \dots, 1; 1, \dots, 1)$, and 
denote $\tilde \mathbb{A} = \mathbb{A}(\mathbf 1, \mathbf 1)$,
$\tilde X^{\pm} = X^{\pm}(\mathbf 1, \mathbf 1)$, 
$\tilde X^0 = X^0(\mathbf 1, \mathbf 1)$ and 
$\tilde Y = Y(\mathbf 1, \mathbf 1)$.
Denote the coordinates of $\tilde \mathbb{A}$ by
$\tilde x_1, \dots, \tilde x_m, \tilde y_1, \dots, \tilde y_n$, and
define $\tilde A_i^{\pm}$ and so on.
Then there is a $G_m$-equivariant finite Galois morphism $\pi_{\mathbb{A}}: 
\tilde \mathbb{A} \to \mathbb{A}$ given by
\[
\tilde x_i \mapsto x_i = \tilde x_i^{a_i}, \quad
\tilde y_j \mapsto y_j = \tilde y_j^{b_j}
\]
with Galois group 
\[
G \cong \prod_i \mathbb{Z}_{a_i} \times \prod_j \mathbb{Z}_{b_j}.
\]
There are induced Galois morphisms
$\pi^{\pm}: \tilde X^{\pm} \to X^{\pm}$, 
$\pi_X^0: \tilde X^0 \to X^0$ and $\pi_Y: \tilde Y \to Y$
with the same Galois group $G$.

$\tilde U_1^- = \tilde X^- \setminus \tilde A_1^-$ is isomorphic to 
$\mathbb{A}^{m+n-1}$ with coordinates
\[
\tilde x_2/\tilde x_1, \dots, \tilde x_m/\tilde x_1, \tilde x_1\tilde y_1, 
\dots, \tilde x_1\tilde y_n.
\]
The quotient $(U_1^-)' = \tilde U_1^-/G'$ for $G' = 
\prod_{i>1} \mathbb{Z}_{a_i} \times \prod_j \mathbb{Z}_{b_j}$
is again isomorphic to 
$\mathbb{A}^{m+n-1}$ with coordinates
\[
(\tilde x_2/\tilde x_1)^{a_2}, \dots, (\tilde x_m/\tilde x_1)^{a_m}, 
(\tilde x_1\tilde y_1)^{b_1}, \dots, (\tilde x_1\tilde y_n)^{b_n}.
\]
Hence the first assertion.

If the action is not small, then there exists an integer $a_1'$ such that
$0 < a_1' < a_1$, $a_1' \vert a_1$, and $a_1'a_2 \equiv \dots \equiv a_1'b_n
\equiv 0 \mod a_1$ except possibly one of the $a_2, \dots, b_n$.
Then $c > 1$, $c_i > 1$ or $c_{m+j} > 1$ for some $i$ or $j$, a contradiction. 
\end{proof}

\begin{Prop}
Let $a = GCD(a_1, \dots, a_m)$ and $b = GCD(b_1, \dots, b_n)$.
Set $a_i = aa_i'$ and $b_j = bb_j'$.
Then the open subset $U_{1,1}$ is isomorphic to the quotient of the 
affine space $\mathbb{A}^{m+n-1}$
by the group $\mathbb{Z}_{a_1'} \times \mathbb{Z}_{b_1'}$ 
whose action is given by the weights
\[
\frac 1{a_1'}(- a_2', \dots, -a_m', b, 0, \dots, 0), \quad
\frac 1{b_1'}(0, \dots, 0, a, - b_2', \dots, -b_n').
\]
Moreover, the action of the group is small in the sense that 
the induced morphism $\mathbb{A}^{m+n-1} \to U_{1,1}$ is etale in codimension 
$1$.
\end{Prop}

\begin{proof}
$\tilde U_{1,1} = \tilde Y \setminus (\tilde A_1 \cup \tilde B_1)$ is 
isomorphic to 
$\mathbb{A}^{m+n-1}$ with coordinates
\[
\tilde x_2/\tilde x_1, \dots, \tilde x_m/\tilde x_1, \tilde x_1\tilde y_1, 
\tilde y_2/\tilde y_1, \dots, \tilde y_n/\tilde y_1.
\]
The quotient $U_{1,1}' = \tilde U_{1,1}/G'$ for $G' = 
\prod_{i>1} \mathbb{Z}_{a_i} \times \prod_{j>1} \mathbb{Z}_{b_j}$
is again isomorphic to 
$\mathbb{A}^{m+n-1}$ with coordinates
\[
(\tilde x_2/\tilde x_1)^{a_2}, \dots, (\tilde x_m/\tilde x_1)^{a_m}, 
\tilde x_1\tilde y_1, 
(\tilde y_2/\tilde y_1)^{b_2}, \dots, (\tilde y_n/\tilde y_1)^{b_n}.
\]
The group $\mathbb{Z}_{a_1} \times \mathbb{Z}_{b_1}$ acts on $U_{1,1}'$ 
with weights
\[
\frac 1{a_1}(- a_2, \dots, -a_m, 1, 0, \dots, 0), \quad
\frac 1{b_1}(0, \dots, 0, 1, - b_2, \dots, -b_n).
\]
The quotient $U_{1,1}'' = U_{1,1}'/G''$ for $G'' = 
\mathbb{Z}_a \times \mathbb{Z}_b$ is still isomorphic to 
$\mathbb{A}^{m+n-1}$ with coordinates
\[
(\tilde x_2/\tilde x_1)^{a_2}, \dots, (\tilde x_m/\tilde x_1)^{a_m}, 
(\tilde x_1\tilde y_1)^{ab}, 
(\tilde y_2/\tilde y_1)^{b_2}, \dots, (\tilde y_n/\tilde y_1)^{b_n}.
\]
Here we note that $GCD(a,b) = 1$.
We can check that there are at least $2$ numbers which are coprime to $a'_1$ 
among $a'_2, \dots, a'_m$ and $b$ by Proposition~\ref{coprime}.
We also check a similar statement for the $b_j'$ and $a$.
Hence the assertion.
\end{proof}

\begin{Defn}
We put natural orbifold structures on the toric varieties 
$X^{\pm}(\mathbf a, \mathbf b)$, and we take the fiber product 
$Y(\mathbf a, \mathbf b)$ in the sense of orbifolds. 
More precisely, the orbifold structure of the latter is not the natural one 
given by the minimal coverings $U_{i,j}'' \to U_{i,j}$ but by
the coverings $U_{i,j}' \to U_{i,j}$.
\end{Defn}

The ideal sheaves 
$\mathcal{O}_{X^-}^{\text{orb}}(- A_i^-)$ and 
$\mathcal{O}_{X^-}^{\text{orb}}(- B_j^-)$
have the structure of invertible orbifold sheaves.
Indeed, on the covering $(U_1^-)'$ of the affine open subset $U_1^-$,
the sheaves $\mathcal{O}_{X^-}^{\text{orb}}(- A_i^-)$ for $i=2,\dots,m$ and 
$\mathcal{O}_{X^-}^{\text{orb}}(- B_j^-)$ for $j=1,\dots,n$ 
are generated by the coordinates
\[
(\tilde x_2/\tilde x_1)^{a_2}, \dots, (\tilde x_m/\tilde x_1)^{a_m}, 
(\tilde x_1\tilde y_1)^{b_1}, \dots, (\tilde x_1\tilde y_n)^{b_n}.
\]
The sheaves of invariants under the Galois group action coincide with the usual 
reflexive ideal sheaves:
\[
\mathcal{O}_{X^-}^{\text{orb}}(- A_i^-)\vert_{X^-}
= \mathcal{O}_{X^-}(- A_i^-), \quad
\mathcal{O}_{X^-}^{\text{orb}}(- B_j^-)\vert{X^-}
= \mathcal{O}_{X^-}(- B_j^-).
\]

We can define invertible orbifold sheaves 
$\mathcal{O}_{X^{\pm}}^{\text{orb}}(k)$ for $k \in \mathbb{Z}$
on $X^{\pm}$ and $\mathcal{O}_Y^{\text{orb}}(k_1,k_2)$
for $k_1,k_2 \in \mathbb{Z}$ on $Y$
so that we have isomorphisms
\[
\begin{split}
&\mathcal{O}_{X^{\pm}}^{\text{orb}}(A_i^{\pm}) \cong 
\mathcal{O}_{X^{\pm}}^{\text{orb}}(\pm a_i), \quad
\mathcal{O}_{X^{\pm}}^{\text{orb}}(B_j^{\pm}) \cong 
\mathcal{O}_{X^{\pm}}^{\text{orb}}(\mp b_j) \\
&\mathcal{O}_Y^{\text{orb}}(A_i) \cong 
\mathcal{O}_Y^{\text{orb}}(a_i, 0), \quad
\mathcal{O}_Y^{\text{orb}}(B_j) \cong 
\mathcal{O}_Y^{\text{orb}}(0, b_j) \\
&\mathcal{O}_Y^{\text{orb}}(\bar E) \cong 
\mathcal{O}_Y^{\text{orb}}(-1, -1) 
\end{split}
\]
where $\bar E$ is the exceptional prime divisor on the Galois covers
$U'_{i,j}$ so that 
\[
\pi_{i,j}^*E = ab\bar E
\]
for $\pi_{i,j}: U'_{i,j} \to U_{i,j}$.
Indeed, the coordinates 
\[
(\tilde x_2/\tilde x_1)^{a_2}, \dots, (\tilde x_m/\tilde x_1)^{a_m}, 
\tilde x_1\tilde y_1, 
(\tilde y_2/\tilde y_1)^{b_2}, \dots, (\tilde y_n/\tilde y_1)^{b_n}.
\]
on $U'_{1,1}$ correspond to the prime divisors $A_2,\dots,A_m$, $\bar E$, and
$B_2,\dots, B_n$.
We have the following equalities
\[
\begin{split}
&(\mu^-)^*A_i^- = A_i, \quad
(\mu^-)^*B_j^- = B_j + b_j\bar E \\
&(\mu^+)^*A_i^+ = A_i + a_i\bar E, \quad
(\mu^+)^*B_j^- = B_j \\
&(\mu^-)^*\mathcal{O}_{X^-}^{\text{orb}}(k) 
= \mathcal{O}_Y^{\text{orb}}(k,0), \quad
(\mu^+)^*\mathcal{O}_{X^+}^{\text{orb}}(k) 
= \mathcal{O}_Y^{\text{orb}}(0,k).
\end{split}
\]
Since $K_{X^{\pm}} + \sum_i A_i^{\pm} + \sum_j B_j^{\pm} 
\sim 0$, we have
\[
\omega_{X^{\pm}}{\text{orb}} \cong 
\mathcal{O}_{X^{\pm}}{\text{orb}}(\pm (\sum_i a_i-\sum_jb_j)).
\]
Hence 
\[
(\mu^-)^*K_{X^-} \sim (\mu^+)^*K_{X^+} + (\sum_i a_i-\sum_jb_j)\bar E.
\]
On the other hand, though we have 
$K_Y + \sum_i A_i + \sum_j B_j + E \sim 0$, we have 
\[
\omega_Y{\text{orb}} \cong 
\mathcal{O}_Y{\text{orb}}(-\sum_i a_i+1,-\sum_jb_j+1)
\]
because we have additional ramification along $E$.


\section{Flip and derived categories}

The following example shows that we should consider the derived categories of 
orbifold sheaves instead of ordinary sheaves.

\begin{Expl}\label{counterexample}
We consider Francia's flop; we take a sequence of integers 
$(\mathbf a; \mathbf b)=(1,2;1,1,1)$.

$X^-$ has only one singular point $P_0 \in U^-_1$ 
which is a quotient
singularity of type $\frac 12(1,1,1,1)$, while $X^+$ is smooth.
$Y$ has $2$-dimensional singular locus of type $A_1$.
The exceptional loci $E^- \subset X^-$, $E^+ \subset X^+$ and
$E \subset Y$ are respectively isomorphic to $\mathbb{P}^1$, 
$\mathbb{P}^2$ and $\mathbb{P}^1 \times \mathbb{P}^2$.

By direct calculation, we observe that the image of a sheaf 
under the Fourier-Mukai transform 
\[
R\mu^+_*L(\mu^-)^*\mathcal{O}_{X^-}(-A_2^-)
\in D^-((X^+)_{\text{coh}})
\] 
has unbounded cohomology sheaves.

On the other hand, we can calculate
\[
R\mu^-_*L(\mu^+)^*(\Omega_{E^+}^1(-1))= 
R\mu^-_*(\mu^+)^!(\Omega_{E^+}^1(1)) = 0
\]
in $D^b((X^-)_{\text{coh}})$.
Indeed, we have
\[
R(\mu^-)_*^{\text{orb}}L(\mu^+)^*_{\text{orb}}(\Omega_{E^+}^1(-1)) = 
R(\mu^-)_*^{\text{orb}}(\mu^+)^!_{\text{orb}}\Omega_{E^+}^1(1) =
\mathcal{O}_{\tilde P_0}^-[-1]
\]
in $D^b((X^-)_{\text{coh}}^{\text{orb}})$,
where $\mathcal{O}_{\tilde P_0}^-$ is the structure sheaf 
$\mathcal{O}_{\tilde P_0}$ of the point $\tilde P_0 \in \tilde U_1^-$ 
above $P_0$
with the non-trivial action of $\text{Gal}(\tilde U_1^-/U_1^-)$.
\end{Expl}

The following theorem is the main result of this section.  
This is an extension of the result of Bondal-Orlov \cite{BO1}~Theorem~3.6 
to the orbifold case.
First we introduce the notation.
Fixing a sequence of positive integers $(\mathbf a; \mathbf b)$, 
we consider the following cartesian diagram of quasi-projective toroidal 
varieties 
\begin{equation}\label{toroidalflip}
\begin{CD}
\mathcal{Y} @>{\hat \mu^+}>> \mathcal{X}^+ \\
@V{\hat \mu^-}VV   @VV{\hat \phi^+}V \\
\mathcal{X}^- @>>{\hat \phi^-}> \mathcal{X}^0
\end{CD}
\end{equation}
whose local models are the product of the 
toric varieties defined in \S4 and a fixed smooth closed subvariety $W$
of $\mathcal{X}^0$:
\begin{equation}\label{productflip}
\begin{CD}
Y(\mathbf a; \mathbf b) \times W @>{\mu^+ \times \text{Id}_W}>> 
X^+(\mathbf a; \mathbf b) \times W \\
@V{\mu^- \times \text{Id}_W}VV   @VV{\phi^+ \times \text{Id}_W}V \\
X^-(\mathbf a; \mathbf b) \times W 
@>>{\phi^- \times \text{Id}_W}> X^0(\mathbf a; \mathbf b) \times W.
\end{CD}
\end{equation}
We assume that the base change of the diagram \ref{toroidalflip} by the
completion of $\mathcal{X}^0$ at any point $w \in W$ 
is isomorphic to that of the diagram \ref{productflip} by the completion of 
$X^0(\mathbf a; \mathbf b) \times W$ at $(P_0, w)$, where 
$P_0 = \phi^{\pm}(E^{\pm})$.
We put natural orbifold structures on $\mathcal{X}^{\pm}$ and 
the orbifld structure of the fiber product on $\mathcal{Y}$.

\begin{Thm}\label{main1}
In the situation above, 
assume that the orbifolds $\mathcal{X}^{\pm}$ have 
Cohen-Macaulay global covers.
Asume moreover that $\sum a_i \le \sum b_j$, i.e., 
\[
(\hat \mu^-)^*K_{\mathcal{X}^-} \le (\hat \mu^+)^*K_{\mathcal{X}^+}.
\]
Then the Fourie-Mukai functors
\[
\begin{split}
&\mathcal{F} = 
R(\hat \mu^+)_*L(\hat \mu^-)^*: D^b_{\text{coh}}(\mathcal{X}^-) 
\to D^b_{\text{coh}}(\mathcal{X}^+) \\
&\mathcal{F}' = 
R(\hat \mu^+)_*(\hat \mu^-)^!: D^b_{\text{coh}}(\mathcal{X}^-) 
\to D^b_{\text{coh}}(\mathcal{X}^+)
\end{split}
\]
are fully faithful.
In particular, if $\sum a_i = \sum b_j$, then they are equivalences 
of categories. 
\end{Thm}

We recall the definition of the spanning class.

\begin{Defn}
A set of objects $\Omega$ of a triangulated category $A$ is
said to be a {\it spanning class} if the following hold 
for any $a \in A$:

(1) $\text{Hom}_A(a, \omega[k]) = 0$ for any $\omega \in \Omega$
and $k \in \mathbb{Z}$ implies $a \cong 0$.

(2) $\text{Hom}_A(\omega[k],a) = 0$ for any $\omega \in \Omega$
and $k \in \mathbb{Z}$ implies $a \cong 0$.
\end{Defn}

\begin{Lem}
Let $f: A \to B$ 
be an exact functor between triangulated categories with a right adjoint $g$ 
and a left adjoint $h$.
Let $\Omega$ be a spanning class of $A$.
Assume that $gf(\omega) \cong hf(\omega) \cong \omega$ 
for any $\omega \in \Omega$.
Then $gf(a) \cong hf(a) \cong a$ for any $a \in A$ and $f$ is fully faithful.
\end{Lem}

\begin{proof}
If $a \in A$, then
\[
\begin{split}
&\text{Hom}_A(\omega, gf(a)) \cong \text{Hom}_A(hf(\omega), a) \cong 
\text{Hom}(\omega, a) \\
&\text{Hom}_A(hf(a), \omega) \cong \text{Hom}_A(a, gf(\omega)) \cong 
\text{Hom}(a, \omega)
\end{split}
\]
hence the natural morphisms $a \to gf(a)$ and $hf(a) \to a$ are isomorphisms. 
Thus
\[
\text{Hom}(f(a), f(a')) \cong \text{Hom}(hf(a), a') \cong 
\text{Hom}(a, a')
\]
for any $a' \in A$.
\end{proof}

\begin{Expl}
(1) Let $X$ be a quasi-projective variety with an orbifold structure
having a Cohen-Macaulay global cover.
For any point $x \in X$,
there exists a finite group $G_x$ such that,
if $x \in \pi_i(X_i)$, then the stabilizer subgroup of $G_i$ at any point 
$\tilde x \in \pi_i^{-1}(x)$ is isomorphic to $G_x$.
Let $V$ be any irreducible representation of $G_x$.
Then the sheaf
\[
Z_{x,V,i} = \bigoplus_{\tilde x \in \pi_i^{-1}(x)} V \otimes_{\mathbb{C}} 
\mathcal{O}_{\tilde x}
\]
on $X_i$ glue together to define an orbifold sheaf $Z_{x,V}$ on $X$.

Let $\mathcal{P}_X$ be the set of all 
the orbifold sheaves of the form $Z_{x,V}$ 
for the points $x \in X$ and irreducible representations $V$ of $G_x$.
Then $\mathcal{P}_X$ is a spanning class for
$D^b(X_{\text{coh}}^{\text{orb}})$.
The proof is similar to \cite{Bridgeland1}~Example~2.2.

(2) Let $X^{\pm} = X^{\pm}(\mathbf a; \mathbf b)$ be as in \S4 and 
fix a positive integer $k_0$. 
Then the set of orbifold sheaves 
\[
\mathcal{Q}_{X^{\pm}} = \{\mathcal{O}_{X^{\pm}}^{\text{orb}}(k)
\vert k \in \mathbb{Z} \text{ and } k \ge k_0\}
\]
is a spanning class for
$D^b((X^{\pm})_{\text{coh}}^{\text{orb}})$.
Indeed, any orbifold sheaf in $\mathcal{P}_{X^{\pm}}$
can be resolved into
a complex of orbifold sheaves which are direct sums 
of the orbifold sheaves in $\mathcal{Q}_{X^{\pm}}$.
\end{Expl}

\begin{Lem}\label{triple}
Let $X^{\pm} = X^{\pm}(\mathbf a; \mathbf b)$ be as in \S4.
\[
\begin{split}
&F = R(\mu^+)_*^{\text{orb}} \circ L(\mu^-)^*_{\text{orb}}:
D^b((X^-)_{\text{coh}}^{\text{orb}}) \to 
D^b((X^+)_{\text{coh}}^{\text{orb}}) \\
&G = R(\mu^-)_*^{\text{orb}} \circ (\otimes \omega_{Y/X^+}^{\text{orb}})
\circ L(\mu^+)^*_{\text{orb}}:
D^b((X^+)_{\text{coh}}^{\text{orb}}) \to 
D^b((X^-)_{\text{coh}}^{\text{orb}}) \\
&H = R(\mu^-)_*^{\text{orb}} \circ (\otimes \omega_{Y/X^-}^{\text{orb}})
\circ L(\mu^+)^*_{\text{orb}}:
D^b((X^+)_{\text{coh}}^{\text{orb}}) \to 
D^b((X^-)_{\text{coh}}^{\text{orb}}) \\
&F' = R(\mu^+)_*^{\text{orb}} \circ (\otimes \omega_{Y/X^-}^{\text{orb}})
\circ L(\mu^-)^*_{\text{orb}}:
D^b((X^-)_{\text{coh}}^{\text{orb}}) \to 
D^b((X^+)_{\text{coh}}^{\text{orb}}) \\
&G' = R(\mu^-)_*^{\text{orb}} \circ 
(\otimes \omega_{Y/X^+}^{\text{orb}}\omega_{Y/X^-}^{\text{orb}-1})
\circ L(\mu^+)^*_{\text{orb}}:
D^b((X^+)_{\text{coh}}^{\text{orb}}) \to 
D^b((X^-)_{\text{coh}}^{\text{orb}}) \\
&H' = R(\mu^-)_*^{\text{orb}} \circ L(\mu^+)^*_{\text{orb}}:
D^b((X^+)_{\text{coh}}^{\text{orb}}) \to 
D^b((X^-)_{\text{coh}}^{\text{orb}}).
\end{split}
\]
Then $(H, F, G)$ and $(H', F', G')$ are adjoint triples of functors.
\end{Lem}

We need a simple lemma in commutative algebra:

\begin{Lem}\label{graded}
Let $R = \mathbb{C}[x_1, \dots, x_m]$ be a polynomial ring
with graded ring structure defined by $\deg x_i =a_i$.
Let $I_k$ be the graded ideal of $R$ consisting of elements of degree
greater than or equal to $k$.
Then there exists a graded free resolution
\[
0 \to \bigoplus_{\lambda \in \Lambda_m} R(-e^{(m)}_{\lambda}) \to \dots \to 
\bigoplus_{\lambda \in \Lambda_1} R(-e^{(1)}_{\lambda}) \to I_k \to 0
\]
given by matrices with monomial entries such that
\[
k \le e^{(l)}_{\lambda} < k + \sum_{i=1}^m a_i
\] 
for any $1 \le l \le m$ and any $\lambda$.
\end{Lem}

\begin{proof}
From the Koszul complex
\[
0 \to R(-\sum_{i=1}^m a_i) \to \dots \to 
\bigoplus_{i=1}^m R(-a_i) \to R \to \mathbb{C} \to 0
\]
we obtain
\[
\text{Ext}^l_R(\mathbb{C}, \mathbb{C}) \cong 
\bigwedge^l (\bigoplus_{i=1}^m \mathbb{C}(a_i)).
\]
We can express the $R$-module $R/I_k$ 
as extensions of the $R$-modules $\mathbb{C}(-e)$
such that $0 \le e < k$.
Since
\[
\text{Ext}^l_R(R/I_k, \mathbb{C}) \cong 
\bigoplus_{\lambda \in \Lambda_l} \mathbb{C}(e^{(l)}_{\lambda})
\]
we obtain our assertion.
\end{proof}

Let $\mathcal{I}_k^-$ ($k \ge 0$) be 
the orbifold ideal sheaf on $X^- = X^-(\mathbf a; \mathbf b)$
generated by monomials of order $k$ on $\tilde y_1, \dots, \tilde y_n$.
By the vanishing theorem, we have the following:

\begin{Lem}
\[
R(\mu^-)_*^{\text{orb}}\mathcal{O}_Y^{\text{orb}}(k\bar E) = 
\mathcal{O}_{X^-}^{\text{orb}}
\]
for $0 \le k \le \sum_j b_j - 1$ and
\[
R(\mu^-)_*^{\text{orb}}\mathcal{O}_Y^{\text{orb}}(-k\bar E) = \mathcal{I}_k^-
\]
for $0 \le k$
\end{Lem}

\begin{Prop}
Under the notation of Lemma~\ref{triple}, 
let $u = \mathcal{O}_{X^-}^{\text{orb}}(k)$.
Assume that $\sum_i a_i \le \sum_j b_j$.

(1) If $k \ge 0$, then
$GF(u) \cong HF(u) \cong u$.

(2) If $k \ge \sum_j b_j-1$, then
$G'F'(u) \cong H'F'(u) \cong u$.
\end{Prop}

\begin{proof}
Since
\[
L(\mu^-)^*_{\text{orb}}\mathcal{O}_{X^-}^{\text{orb}}(k)
\cong \mathcal{O}_Y^{\text{orb}}(0,-k)(-k \bar E)
\]
we have
\[
F(\mathcal{O}_{X^-}^{\text{orb}}(k))
\cong I^+_k(-k).
\]
We have a locally free resolution
\[
0 \to \bigoplus_{\lambda \in \Lambda_m} 
\mathcal{O}_{X^+}^{\text{orb}}(e^{(m)}_{\lambda}) \to 
\dots \to 
\bigoplus_{\lambda \in \Lambda_1} 
\mathcal{O}_{X^+}^{\text{orb}}(e^{(1)}_{\lambda}) \to I^+_k \to 0
\]
where the maps are given by matrices whose entries are 
monomials in the $x_i$.
Therefore
\[
\begin{split}
& (\otimes \omega_{Y/X^+}^{\text{orb}}) \circ
L(\mu^+)^*_{\text{orb}} \circ F(\mathcal{O}_{X^-}^{\text{orb}}(k)) \\
&\cong 
(0 \to \bigoplus_{\lambda \in \Lambda_m} 
\mathcal{O}_Y^{\text{orb}}(k-e^{(m)}_{\lambda},0)
((k-e^{(m)}_{\lambda}+\sum_i a_i-1)\bar E) \to \\
&\dots \to 
\bigoplus_{\lambda \in \Lambda_1} 
\mathcal{O}_Y^{\text{orb}}(k-e^{(1)}_{\lambda},0)
((k-e^{(1)}_{\lambda}+\sum_i a_i-1)\bar E) \to 0).
\end{split}
\]
Since $\sum_i a_i \le \sum_j b_j$ and 
$k \le e^{(l)}_{\lambda} < k + \sum_i a_i$ 
for any $1 \le l \le m$ and any $\lambda$, 
we obtain
\[
\begin{split}
&G \circ F(\mathcal{O}_{X^-}^{\text{orb}}(k)) \\
&\cong 
(0 \to \bigoplus_{\lambda \in \Lambda_m} 
\mathcal{O}_{X^-}^{\text{orb}}(k-e^{(m)}_{\lambda}) \to 
\dots \to 
\bigoplus_{\lambda \in \Lambda_1} 
\mathcal{O}_{X^-}^{\text{orb}}(k-e^{(1)}_{\lambda}) \to 0) \\
&\cong 
\mathcal{O}_{X^-}^{\text{orb}}(k)
\end{split}
\]
where the latter isomorphism is obtained because 
$I_k$ in Lemma~\ref{graded} is primary to the maximal ideal $I_1$.

Other isomorphisms are proved similarly.
\end{proof}

\begin{Cor}\label{key}
If $\sum_i a_i \le \sum_j b_j$, then 
\[
GF(u) \cong HF(u) \cong G'F'(u) \cong H'F'(u) \cong u 
\]
for any $u \in D^b((X^-)_{\text{coh}}^{\text{orb}})$, and 
the functors $F$ and $G$ are fully faithful.
\end{Cor}

\begin{Cor}
If $\sum_i a_i = \sum_j b_j$, then the functor $F$
is an equivalence of categories whose inverse is given by $G$.
\end{Cor}

\begin{proof}[Proof of Theorem~\ref{main1}]
We can extend the result of Corollary~\ref{key} by replacing
$X^{\pm}$ and $Y$ by $X^{\pm} \times W$ and $Y \times W$.
If we define $\mathcal{F}$, $\mathcal{G}$, and so on as in Lemma~\ref{triple},
then we have 
\[
\mathcal{G}\mathcal{F}(u) \cong \mathcal{H}\mathcal{F}(u) 
\cong \mathcal{G}'\mathcal{F}'(u) \cong \mathcal{H}'\mathcal{F}'(u) 
\cong u 
\]
for any $u \in \mathcal{P}_{\mathcal{X}^-}$, hence the result.
\end{proof}


\section{Reconstruction}

We extend the reconstruction theorem by Bondal-Orlov to the orbifold case 
in this section.

\begin{Thm}\label{main2}
Let $X$ and $X'$ be projective varieties with only quotient singularities.
Assume the following conditions for $X$ and $X'$:

(1) The natural orbifold structure has a Cohen-Macaulay global cover.

(2) The canonical divisor generates local class groups at any 
point.

Suppose that $K_X$ or $-K_X$ is ample, and
there is an equivalence of categories 
$D^b(X_{\text{coh}}^{\text{orb}}) \to D^b((X')_{\text{coh}}^{\text{orb}})$
which is compatible with shifting functors.
Then there exists an isomorphism $X \to X'$.
\end{Thm}

\begin{proof}
We follow closely the proof of \cite{BO1}~Theorem~4.5.
Denoting $X$ or $X'$ by $Y$, we let $D(Y) = D^b(Y_{\text{coh}}^{\text{orb}})$ 
and $S_Y$ its Serre functor.

{\it Step 1}. 
We define a {\it point object of codimension $s$} to be an object $P \in D(Y)$ 
satisfying the following conditions:

(1) $S_Y^r(P) \cong P[rs]$ for some positive integer $r$.
Let $r_P$ be the smallest such $r$.

(2) $\text{Hom}^{<0}(P, S_Y^m(P)[-ms]) = 0$ for any integer $m$.

(3) $\text{Hom}^0(P, P) = \mathbb{C}$ and
$\text{Hom}^0(P, S_Y^m(P)[-ms]) = 0$ for $0 < m < r_P$.

{\it Step 1.1}. 
We claim that any point object $P$ of $D(X)$ is of the form 
$\mathcal{O}_x \otimes (\omega_X^{\text{orb}})^m[t]$ 
for some $x \in X$ and $m,t \in \mathbb{Z}$.

Indeed, it follows from (1) that $H^i(P) \otimes 
(\omega_X^{\text{orb}})^r \cong H^i(P)$ and $s = \dim X$, 
hence the supports of the cohomology orbifold sheaves $H^i(P)$ are
$0$-dimensional and $\omega_X^r$ is invertible there.
Then $P$ can be represented by a complex of orbifold sheaves 
whose supports are 
$0$-dimensional.
By (3), the support of $P$ is a single point.  
Let $i_0$ and $i_1$ be the minimum and the maximum of the $i$ such that
$H^i(P) \ne 0$.
Then we have $\text{Hom}^{i_0-i_1}(P, P \otimes 
(\omega_X^{\text{orb}})^m) \ne 0$ for some $m$, hence 
$i_0=i_1$ by (2), i.e., $P$ is an orbifold sheaf.
If the length of $P$ is more than $1$, then there is a non-invertible 
homomorphism
$P \to P \otimes (\omega_X^{\text{orb}})^m$ for some $m$,
a contradiction to (3), hence $P$ has the claimed form.
 
{\it Step 1.2}. 
We claim that any point object $P$ of $D(X')$ is also of the form 
$\mathcal{O}_{x'}\otimes (\omega_{X'}^{\text{orb}})^m[t]$ 
for some $x' \in X'$ and $m,t \in \mathbb{Z}$.

Indeed, since $D(X)$ and $D(X')$ are equivalent, 
it follows from Step 1.1 that, for any point objects $P$ and $Q$ of $D(X')$, 
either 
$Q = S_{X'}^m(P)[t]$ for some $m,t \in \mathbb{Z}$ 
or $\text{Hom}^i(P,Q) = 0$ for any $i$ holds. 
If $P$ is not of the claimed form, 
then $\text{Hom}^i(P, \mathcal{O}_{x'}\otimes (\omega_{X'}^{\text{orb}})^m) 
= 0$ for any $i$, $m$ and $x'$, hence $P = 0$.  

{\it Step 2}. 
An {\it invertible object} $L \in D(Y)$ is defined by the following condition:
If $P$ is any point object of codimension $s$, 
then there exist uniquely integers $m_0 \in [0,r_P-1]$ and $t_0$ 
such that $\text{Hom}^{t_0}(L, S_Y^{m_0}(P)) = \mathbb{C}$ and
$\text{Hom}^i(L, S_Y^m(P)) = 0$ for $i \ne t_0$ or $m \not\equiv m_0 \mod r_P$.

We claim that any invertible object of $D(Y)$ is of the form
$L[t]$ for some invertible orbifold sheaf $L$ on $Y$ and some 
$t \in \mathbb{Z}$. 
Indeed, we consider a spectral sequence
\[
E_2^{p,q} = \text{Hom}^p(H^{-q}(L), \mathcal{O}_y \otimes 
(\omega_Y^{\text{orb}})^m)
\Rightarrow \text{Hom}^{p+q}(L, \mathcal{O}_y \otimes 
(\omega_Y^{\text{orb}})^m)
\]
for $y \in Y$ and $m \in \mathbb{Z}$.
If $i_1$ is the maximum of the $i$ such that
$H^i(L) \ne 0$, then $E_2^{p,-i_1}$ for $p=0,1$ survive at $E_{\infty}$.
On the other hand, for any $y \in \text{Supp}(H^{i_1}(L))$,
there exists an integer $m_y$ such that
$\text{Hom}^0(H^{i_1}(L), \mathcal{O}_y \otimes 
(\omega_Y^{\text{orb}})^{m_y}) \ne 0$.
Thus
\[
\begin{split}
&\text{Hom}^0(H^{i_1}(L), \mathcal{O}_y \otimes 
(\omega_Y^{\text{orb}})^m) = \begin{cases} \mathbb{C} &\text{ if } m \equiv m_y
\mod r_y \\ 
0 &\text{ otherwise} \end{cases} \\
&\text{Hom}^1(H^{i_1}(L), \mathcal{O}_y \otimes 
(\omega_Y^{\text{orb}})^m) = 0
\end{split}
\]
for any $m$, hence $H^{i_1}(L)$ is an invertible orbifold sheaf.
Then $E_2^{p,-i_1} = 0$ for $p \ne 0$, hence $E_2^{0,-i_1+1}$ survives at 
$E_{\infty}$.
Thus $\text{Hom}^0(H^{i_1-1}(L), \mathcal{O}_y \otimes 
(\omega_Y^{\text{orb}})^m) = 0$ for any $y \in Y$ and $m$, 
and $H^{i_1-1}(L) = 0$.
Continuing this process, we obtain that $H^i(L) = 0$ for $i \ne i_1$,
and conclude that $L[i_1]$ is an invertible orbifold sheaf.

{\it Step 3}. 
We fix an invertible orbifold sheaf $L_0$ on $X$.
By Step 2, there exists an invertible orbifold sheaf $L'_0$ on $X'$
such that $L_0 \in D(X)$ corresponds to $L'_0[t_0] \in D(X')$ for some $t_0$.
If we compose the shift functor to the given 
equivalence functor $D(X) \to D(X')$, we may assume that $t_0=0$. 
The set of point objects $P \in D(X)$ such that 
$\text{Hom}(L_0, P) \cong \mathbb{C}$ corresponds bijectively to 
the set of those $P' \in D(X')$ such that $\text{Hom}(L'_0, P') \cong \mathbb{C}$.
They correspond respectively to the sets of points $x \in X$ and $x' \in X'$ 
by the isomorphisms
$P \cong P_x = \mathcal{O}_x \otimes (\omega_X^{\text{orb}})^{m_x}$ and 
$P' \cong P_{x'} = \mathcal{O}_{x'} 
\otimes (\omega_{X'}^{\text{orb}})^{m_{x'}}$,
where $m_x$ and $m_{x'}$ are integers depending on the points $x$ and $x'$. 
Therefore we obtain a bijection of sets of points on $X$ and $X'$.
From this it follows that the sets of invertible orbifold sheaves on 
$X$ and $X'$ also correspond bijectively.

{\it Step 4}.
If $L_1$ and $L_2$ are invertible orbifold sheaves on $X$ 
and $u \in \text{Hom}(L_1,L_2)$, then
the set $U(L_1,L_2,u)$ of points $x \in X$ such that the map  
$u^*: \text{Hom}(L_2 \otimes (\omega_X^{\text{orb}})^m,  P_x) 
\to \text{Hom}(L_1 \otimes (\omega_X^{\text{orb}})^m,  P_x)$ is bijective 
for any $m$ is an affine open subset of $X$.
The subsets $U(L_1,L_2,u)$ for all the $L_1, L_2$ and $u$ 
form a basis of the Zariski topology of $X$.
Hence the Zariski topologies on $X$ and $X'$ coincide under the bijection given
in Step 3.

{\it Step 5}.
Let $L_m = L_0 \otimes  (\omega_X^{\text{orb}})^m$ and
$L'_m = L'_0 \otimes  (\omega_{X'}^{\text{orb}})^m$.
We set $\epsilon = \pm 1$ such that $\epsilon K_X$ is ample. 
Then the subsets $U(L_0,L_{m\epsilon}, u)$ for all the positive integers $m$ 
and all the non-zero sections $u \in \text{Hom}(L_0, L_m)$ 
form a basis of the Zariski topology of $X$.
Since the same statement holds for $X'$, it follows that $\epsilon K_{X'}$ is also
ample. Since $\text{Hom}(L_i, L_{i+m}) \cong H^0(X, mK_X)$ for any $i$,
the multiplication on the (anti-)canonical ring 
$R(X) = \bigoplus_{m=0}^{\infty}H^0(X, m\epsilon K_X)$ 
is given by the composition 
of morphisms in $D(X)$.  Hence $X$ and $X'$ have isomorphic (anti-)canonical 
rings.
\end{proof}


Department of Mathematical Sciences, University of Tokyo, 

Komaba, Meguro, Tokyo, 153-8914, Japan 

kawamata@ms.u-tokyo.ac.jp

\end{document}